\documentclass[12pt]{amsart}
\usepackage{amssymb,amsfonts,amsmath,amsopn,url,bbm,amscd}
\newtheorem{thm}{Theorem}[section]
\newtheorem{prop}[thm]{Proposition}
\newtheorem{lemma}[thm]{Lemma}

\theoremstyle{definition}
\newtheorem*{notation}{Notation}
\newtheorem*{definition}{Definition}
\newcommand{\field}[1]{\mathbbm{#1}}

\newcommand{\Z}{\field{Z}}

\newcommand{\ideal}[1]{\mathfrak{#1}}
\newcommand{\m}{\ideal{m}}
\newcommand{\n}{\ideal{n}}

\newcommand{\intcl}[1]{{#1}^{-}}
\newcommand{\sptc}[1]{{#1}^{*\rm sp}}

\title
{A tight closure analogue of analytic spread}

\author{Neil M. Epstein}

\address{Department of Mathematics,
University of Kansas, Lawrence, KS 66045}

\email{epstein@math.ukans.edu}

\date{June 1, 2004}

\begin{document}

\begin{abstract}
An analogue of the theory of integral closure and reductions is developed for
a more general class of closures, called Nakayama closures.  It is shown
that tight closure is a Nakayama closure by proving a ``Nakayama lemma for tight closure''.
Then, after strengthening A. Vraciu's theory of $*$-independence and the special part of tight closure,
it is shown that all minimal $*$-reductions of an ideal in an analytically irreducible
excellent local ring of positive characteristic have the same minimal number of generators.
This number is called the $*$-spread of the ideal, by analogy with the notion of analytic spread.
\end{abstract}
\maketitle

\section{Introduction}
Fifty years ago, Northcott \& Rees \cite{NR} developed a theory of integral closure and reductions
of ideals.  In particular, they proved that any ideal has minimal reductions, and that if the residue
field is infinite,
all minimal reductions of an ideal have
the same size minimal generating sets.  This common
number is called the ``analytic spread'' of the ideal.  Along the way they proved
that the reduction of an ideal in such a ring is minimal if and only if it has an analytically independent generating
set. In this paper, we prove analogous results for tight closure theory.

Melvin Hochster and Craig Huneke's theory of tight closure has proved extremely useful and powerful, especially
when dealing with a Noetherian ring that contains a field.  As a basic
reference for unexplained notation or terminology, see Huneke's monograph \cite{HuTC}.
All the rings in this paper are Noetherian, local, and of prime characteristic $p>0$.
Given an ideal $I$, an element $x \in R$ is in the \emph{tight closure} of $I$
(and we write $x \in I^*$), if there exists some power $q_0$ of $p$ and some $c$ not in any minimal
prime of $R$, \emph{both possibly dependent
on $x$ and $I$}, such that for all powers $q \geq q_0$ of $p$, $c x^q \in I^{[q]}$.  We say that an
element $c$ is a \emph{$q_0$-weak test element} (or just a \emph{weak test element} for short)
of $R$ if it is not in any minimal prime of $R$ and if for \emph{all} $x$ and \emph{all} $I$, $x \in I^*$ if and
only if for all powers $q \geq q_0$ of $p$, $c x^q \in I^{[q]}$.  In
\cite[Theorem 6.1]{HHbase}, Hochster and Huneke prove the highly non-obvious fact
that every excellent Noetherian local ring of characteristic $p>0$ has a weak test element.  Moreover,
any reduced such ring has a \emph{test element} (i.e. a $1$-weak test element).  We shall use these
facts repeatedly without further comment.

To better understand both tight and integral closure, it is natural to compare the two theories
to see which parts can be applied to one another.
In Section~\ref{sec:*red}, we show how the methods of Northcott and Rees may be used to prove
both the existence of ``minimal $*$-reductions'' (and more generally of ``minimal $cl$-reductions'' for
a wide variety of closure operations $cl$) and their close relationship with the notion
of $*$-independence, a partial analogue of analytic independence for tight closure theory introduced by Adela Vraciu \cite{Vr*ind}.
In order to do so, we prove a ``Nakayama lemma for tight closure theory,'' a tool
which we use in other work in progress and believe will prove very useful to tight closure theorists in the future.

The main theorem of this note, Theorem~\ref{thm:*spread},
states that ``$*$-spread'' (an analogue to analytic spread in tight closure theory) is well-defined in
excellent analytically irreducible local rings of positive characteristic.  That is, the minimal number
of generators required to generate an arbitrary minimal $*$-reduction of an ideal $I$ in such a ring
is an invariant, (called the $*$-spread) of $I$.

In order to show this, we
use an idea originated by Vraciu in \cite{Vr*ind} and further developed by Huneke and Vraciu in \cite{HuVr}.  In \cite{Vr*ind},
Vraciu distinguishes a ``special'' part of the tight closure of an ideal.  She shows that in
certain rings, the tight closure may be represented as a sum of
the ideal with this ``special'' part of its tight closure.  In \cite{HuVr},
Huneke and Vraciu extend Vraciu's result to a larger class
of rings.  In this paper, we extend it to a still larger class of rings and reinterpret its
significance to make it more powerful as a tool.  Namely, it seems enough, for this to
be a useful tool, that the tight closure of an ideal be a kind of ``$q$'th root'' of such a sum.  We
prove in Theorem~\ref{thm:huvrq} that this situation occurs for excellent analytically irreducible local rings of positive characteristic which have the property that the residue field of the normalization coincides with the
original residue field.  (In~\cite{nmeVr},
Vraciu and the present author will use this theorem to give a characterization of $*$-spread in terms of length.)

Also, we use the following proposition of Ian Aberbach frequently in this paper.  It affords a
great deal of control over what happens to elements \emph{not} in the tight closure of an ideal.

\begin{prop}[``colon criterion'']\cite[Proposition 2.4]{AbflatFreg}\label{prop::}
Let $(R,\m)$ be an excellent, analytically irreducible Noetherian local ring of characteristic $p>0$,
let $I$ be an ideal, and let $f \in R$.  Then if $f \notin I^*$, there exists $q_0 = p^{e_0}$
such that for all powers $q \geq q_0$ of $p$, we have \[
I^{[q]} : f^q \subseteq \m^{[q / q_0]}.
\]
\end{prop}

\section{$*$-reductions and a Nakayama lemma for tight closure theory}\label{sec:*red}
Let $(R,\m)$ be a Noetherian local ring.  In general,
let $cl$ be a closure operation on ideals, in the sense that
for any ideal $I$, we have that $I^{cl}$ is an ideal, $I \subseteq I^{cl} = (I^{cl})^{cl}$,
and for any ideals $I$ and $J$, if $J \subseteq I$ then $J^{cl} \subseteq I^{cl}$.
Clearly integral closure is a closure in this sense, and if $R$ is of equal characteristic,
so is tight closure.

We say that $cl$ is a \emph{Nakayama closure} if for any ideals $I$ and $J$
such that $J \subseteq I \subseteq (J + \m I)^{cl}$, it follows that $I \subseteq J^{cl}$.

It is easy to see that integral closure is Nakayama \cite{NR}.  For suppose that \[
J \subseteq I \subseteq \intcl{(J + \m I)}.
\]
Then there is some integer $r$ such that \[
I^{r+1} = (J + \m I) I^r = J I^r + \m I^{r+1}.
\]
Then the Nakayama lemma implies that $I^{r+1} = J I^r$, which in turn implies
that $I \subseteq \intcl{J}$.

\begin{notation}
For an ideal $I$, the symbol $\intcl{I}$ will be used to denote the integral closure of the ideal.
\end{notation}

It is also true, under reasonable hypotheses, that tight closure is Nakayama, as we show
in the following proposition.\footnote{For a more powerful
version (and application) of this proposition, see \cite{nmepdep}.}

\begin{prop}[Nakayama lemma for tight closure]\label{prop:Nak*}
Let $(R,\m,k)$ be a Noetherian local ring of characteristic $p>0$ possessing
a $q_0$-weak test element $c$.  Let
$I$, $J$ be ideals of $R$ such that
$J \subseteq I \subseteq (J + \m I)^*$.  Then $I \subseteq J^*$.
\end{prop}

\begin{proof}
Let $c_1 = c$, and for all integers $r \geq 1$, inductively define $c_{r+1} = c c_r^{q_0}$.
We begin by proving the following

\noindent \textbf{Claim: } For any integer $r \geq 1$, \begin{equation}\label{eq:qrclaim}
c_r I^{[q]} \subseteq J^{[q]} + (\m^{[q]})^r I^{[q]}
\end{equation}
for all powers $q \geq q_0^r$ of $p$.

\begin{proof}[Proof of Claim.]
The case where $r=1$ follows from the hypothesis along with the definition of a $q_0$-weak test
element.  So let $r \geq 1$, assume we have proved (\ref{eq:qrclaim}) for $r$, and let $q \geq q_0^{r+1}$ be
a power of $p$.  Let $q' = q / q_0$, and note that $q' \geq q_0^r$.  Then: \begin{align*}
c_{r+1} I^{[q]} &= c (c_r I^{[q']})^{[q_0]} \\
&\subseteq c \left(J^{[q']}\right)^{[q_0]} + c \left((\m^{[q']})^r I^{[q']}\right)^{[q_0]} \\
&= c J^{[q]} + \left(\m^{[q]}\right)^r \left(c I^{[q]}\right) \\
&\subseteq c J^{[q]} + \left(\m^{[q]}\right)^r \left(J^{[q]} + \m^{[q]} I^{[q]}\right) \\
&\subseteq J^{[q]} + \left(\m^{[q]}\right)^{r+1} I^{[q]}.
\end{align*}
\end{proof}
For a fixed $r$, the fact that (\ref{eq:qrclaim}) holds for all $q \geq q_0^r$ implies
that $I \subseteq (J + \m^r I)^*$.  Since this holds for all $r \geq 1$, we have that
for all $q \geq q_0$, \[
c I^{[q]} \subseteq \bigcap_{r \geq 1} \left( J^{[q]} + \left(\m^{[q]}\right)^r I^{[q]} \right) \subseteq J^{[q]},
\]
where the final containment follows from the Krull intersection theorem.  Thus, $I \subseteq J^*$.
\end{proof}

For a closure operation $cl$, we say that $J$ is a \emph{$cl$-reduction} of $I$
if $J \subseteq I \subseteq J^{cl}$.  We say it is \emph{minimal} if it is minimal
with respect to inclusion.

The proofs of several statements in \cite{NR} about integral closure are in fact
valid for any Nakayama closure.  Most notably:
\begin{lemma}\label{lem:mingextend}
If $cl$ is a Nakayama closure, then for any $cl$-reduction $J$ of $I$, there
is a minimal $cl$-reduction $K$ of $I$ contained in $J$.  Moreover, in this situation
any minimal generating set of $K$ extends to a minimal generating set of $J$.
\end{lemma}

\begin{proof}
(Essentially in \cite{NR}.)  $V = I / \m I$ is a finite-dimensional $k$-vector space.  Consider the set
of ideals \[
\Sigma = \{L + \m I \mid L \subseteq J \text{ and } I \subseteq (L + \m I)^{cl} \}.
\]
For any $L' \in \Sigma$, $L' / \m I$ is a subspace of $V$.  Moreover, $\Sigma$
is nonempty: It contains $J+ \m I$, since $I \subseteq J^{cl} \subseteq (J + \m I)^{cl}$.
Choose an ideal $L \subseteq J$ among such ideals
so that the corresponding subspace $W = (L + \m I) / \m I$ of $V$ has the smallest possible dimension, say $\ell$.  Then
choose $f_1, \dotsc, f_\ell \in L$ such that their images in $W$ form a basis for $W$.  Let $K = (f_1, \dotsc, f_\ell)$.  Note first that $K$ is a $cl$-reduction of $I$, since we have \[
K + \m I \subseteq I \subseteq (L + \m I)^{cl} = (K + \m I)^{cl},
\]
and then the fact that $cl$ is Nakayama implies that $I \subseteq K^{cl}$.  As for minimality, suppose that
$T \subseteq K$ is an ideal such that $I \subseteq T^{cl}$.  Then $X = (T + \m I) / \m I$ is a subspace of $W$, so
its dimension must be the same since $W$ was chosen with the smallest possible dimension.  Therefore $W = X$.
So $T + \m I = K + \m I$, so that $K \subseteq T + \m I$.  But then for any $k \in K$, we have $k = t + x$ for
some $t \in T$ and $x \in \m I$, in which case $x = k-t \in K \cap \m I = \m K$.  The last equality follows from
the fact that the images of $f_1, \dotsc, f_\ell$ are a basis for $W$.  Thus, $K \subseteq T + \m K$, so
that the standard Nakayama lemma implies that $K = T$, whence $K$ is a minimal $cl$-reduction of $I$.

For the second statement of the lemma, let $K$ be a minimal $cl$-reduction of $I$
contained in $J$.  Let $L \subseteq K$ such that $L + \m I = K + \m I$.  Then $I \subseteq (K + \m I)^{cl}
= (L + \m I)^{cl}$, so $I \subseteq L^{cl}$ since $cl$ is a Nakayama closure.  Then minimality of
$K$ with this property forces $L = K$.  Hence, any set of elements of $K$ whose images form a basis of
$W = (K + \m I) / \m I$ must generate all of $K$, and indeed they must form a minimal set of generators of $K$.
Hence, $K / \m K \cong (K + \m I) / \m I \cong K / (K \cap \m I)$, so we have $K \cap \m I = \m K$.
Thus also $K \cap \m J = \m K$, which means that any minimal set of generators of $K$ extends to one of $J$.
\end{proof}

We say (mimicking Vraciu's definition of $*$-independence in~\cite{Vr*ind})
that a set of elements $x_1, \dotsc, x_n \in R$ is \emph{$cl$-independent} if
for $1 \leq i \leq n$, $x_i \notin (x_1, \dotsc, \hat{x_i}, \dotsc, x_n)^{cl}$.
We say that an ideal is $cl$-independent if it has a $cl$-independent
generating set.  It is obvious that any such set must be a minimal generating
set for the ideal it generates.  Accordingly, we say that an ideal is
\emph{strongly} $cl$-independent if \emph{every} minimal set of generators of the ideal
is $cl$-independent.

In the case of integral closure, independence and strong independence can
differ, even for normal hypersurfaces.  For example, let $S = k[[X,Y,Z,W]]$
be the power series ring in 4 variables over
a field $k$, let $R = S / (X Y - Z W)$, and let the lower-case letters $x,y,z,w$ denote elements of $R$.
Then the set $x,y,z,w$ in the ring $R$ is
$^-$-independent but not strongly $^-$-independent.  To see that they are $^-$-independent,
we need only note that the ideals $(x,y,z)$, $(x,y,w)$, $(x,z,w)$ and $(y,z,w)$ are prime,
and hence integrally closed.
However, the set
$x+y, y, z, w$, which generates the same ideal as $x,y,z,w$, is not $^-$-independent,
since $y^2 = (x+y) y - z w$ shows that $y \in \intcl{(x+y,z,w)}$.

On the other hand, Vraciu shows in~\cite[Proposition 3.3]{Vr*ind}
that for ideals in excellent analytically irreducible
local rings of characteristic $p>0$, $*$-independence and strong $*$-independence coincide.

Perhaps surprisingly, the notions of strong independence and minimal reduction are
intimately related, according to the following result.  Thanks to Dan Katz for
the idea of how to prove this proposition in the case of integral closure.

\begin{prop}\label{prop:strongind}
Let $cl$ be a Nakayama closure, and let $J$ be a $cl$-reduction of $I$.
Then $J$ is a minimal $cl$-reduction of $I$ if and only if $J$ is a strongly $cl$-independent ideal.
\end{prop}

\begin{proof}
Let $z_1, \dotsc, z_n$ be a minimal set of generators of $J$.  If, say, $z_1 \in (z_2, \dotsc, z_n)^{cl}$, then $(z_2, \dotsc, z_n)$ is also a $cl$-reduction
of $I$, and strictly smaller than $J$ by the minimality of the set $z_1, \dotsc, z_n$ of generators.  Hence,
if $J$ is a \emph{minimal} $cl$-reduction of $I$, then any minimal set of generators $z_1, \dotsc, z_n$ of $J$ is $cl$-independent.  That is, $J$ is a strongly $cl$-independent ideal.

Conversely, suppose that $J$ is a strongly $cl$-independent ideal and $K \subseteq J$ is
a minimal $cl$-reduction of $I$.  By Lemma~\ref{lem:mingextend}, any minimal
set of generators of $K$ may be extended to form a minimal set of generators of $J$.
That is, there are elements $y_1, \dotsc, y_n$ and some $r \leq n$ such that $y_1, \dotsc, y_r$
form a minimal generating set for $K$, while $y_1, \dotsc, y_n$ form a minimal generating set
for $J$.  Since $y_1, \dotsc, y_n$ are $cl$-independent, $y_n \notin (y_1, \dotsc, y_{n-1})^{cl}$, so that
if $n > r$, $y_n \notin (y_1, \dotsc, y_r)^{cl} = K^{cl} = J^{cl}$.  But $y_n \in J \subseteq J^{cl}$, which
is a contradiction.  Therefore $n=r$, so that $J=K$ is a minimal $cl$-reduction of $I$.
\end{proof}

By analogy with the notion of analytic spread, we make the following

\begin{definition}
Let $cl$ be a Nakayama closure operation.  We say that an ideal $I$ of $R$ \emph{has $cl$-spread}
if all minimal $cl$-reductions of $I$ have the same size minimal
generating sets.  We call this common number the \emph{$cl$-spread of $I$}, denoted $\ell^{cl}(I)$.
\end{definition}

In particular, Northcott and Rees showed that if $R$ has infinite residue field then the analytic
spread of an ideal is its $^-$-spread.
However, I don't know in general
whether ideals of Noetherian local rings with finite residue field have $^-$-spread or not.\footnote{
Zariski and Samuel give an example in an
appendix to \cite{ZS2} of a Noetherian local ring $R$ with finite residue field and an ideal $I$ with
analytic spread 1 but no principal reductions.  However, their ideal $I$ has no
proper reductions at all and is minimally generated by 2 elements, so it has $^-$-spread 2.}

\begin{prop}
Let $cl$ be a Nakayama closure on ideals of $R$.  If an ideal $I$ of $R$ has $cl$-spread,
then a $cl$-reduction $J$ of $I$ is minimal
if and only if $\mu(J) \leq \ell^{cl}(I)$.
\end{prop}

\begin{proof}
Let $\ell = \ell^{cl}(I)$.  By assumption any minimal $cl$-reduction is minimally generated
by $\ell$ elements.  Conversely, if a $cl$-reduction
$J$ of $I$ can be generated by $\ell$ or fewer elements, let $K$ be a minimal $cl$-reduction
of $I$ inside of $J$.  Then $K$ is minimally generated by $\ell$ elements $x_1, \dotsc, x_\ell$,
which can be extended to form a minimal generating set for $J$.  But since $\mu(J) \leq \ell$, this
forces $J=K$.
\end{proof}

\section{The special part of the tight closure}

To find conditions under which all minimal $*$-reductions of an ideal have the same
minimal number of generators (i.e. conditions under which $*$-spread holds), we need
the following tool:

\begin{definition}[Vraciu \cite{Vr*ind}]
For an ideal $I$ in a Noetherian local ring $R$, we define the \emph{special part of the tight closure}\footnote{In Vraciu \cite{Vr*ind} and Huneke \& Vraciu \cite{HuVr}, this
is called the ``special tight closure'' of $I$.  However, since it is not in fact a closure
(c.f. the next lemma),
we prefer the present term.}\ \ of $I$ to be the set \[
\sptc{I} = \{z \in R \mid \exists q \text{ such that } z^{q} \in
(\m I^{[q]})^* \}
\]
\end{definition}

Let us note some general facts about $*$-independence and the special part of tight closure in addition to those
stated in Vraciu \cite{Vr*ind}.  For example,

\begin{lemma}
Let $(R,\m)$ be a Noetherian local ring of characteristic $p>0$ containing a weak
test element $c$.  Let $I$ be an ideal of $R$.  If $I \subseteq \sptc{I}$, then
$I$ is nilpotent.
\end{lemma}

\begin{proof}
Let $f_1, \dotsc, f_n$ be a generating set for $I$.  Then there is some $q$ such
that for $1 \leq i \leq n$, $f_i^{q} \in (\m I^{[q]})^*$.  Hence, $I^{[q]}
\subseteq (\m I^{[q]})^*$.  Then by Proposition~\ref{prop:Nak*}, $I^{[q]} \subseteq 0^*$.
But the tight closure of the zero ideal is the nilradical of the ring.  Hence, every
element of $I^{[q_0]}$ (and thus every element of $I$) is nilpotent.
\end{proof}

Next, we prove a general result, which Vraciu proves in the case where $R$ is excellent and analytically irreducible.
\begin{lemma}\label{lem:sp*capI}
If $R$ is a Noetherian local ring with a weak test element $c$, then
for any ideal $I$ of $R$, $\m I \subseteq \sptc{I} \cap I$, and the reverse containment
holds if $I$ is a $*$-independent ideal.
\end{lemma}

\begin{proof}
The first containment is obvious.  Next, suppose that $I$ is $*$-independent, let
$z_1, \dotsc, z_n$ be $*$-independent generating
set for $I$, and let $z \in \sptc{I} \cap I$.  Then
$z = r_1 z_1 + \cdots r_n z_n$ for some $r_i \in R$.  If $z \notin \m I$, then without
loss of generality $r_1 \notin \m$.  Let $J = (z_2, \dotsc, z_n) \subseteq I$.

By definition of the special part of the tight closure, there is some power $q$ of $p$ such
that $z^{q} \in (\m I^{[q]})^*$.  That is, \[
r_1^{q} z_1^{q} + \sum_{i=2}^n r_i^{q} z_i^{q} \in (\m I^{[q]})^*.
\]
Dividing both sides by the unit $r_1^{q}$ and moving the $z_i$'s other than $z_1$
to the right hand side, we have \[
z_1^{q} \in (z_2^{q}, \dotsc, z_n^{q}) + (\m I^{[q]})^* = J^{[q]} + (\m I^{[q]})^*.
\]
Hence, \[
I^{[q]} \subseteq J^{[q]} + (\m I^{[q]})^*
\subseteq \left( J^{[q]} + \m I^{[q]} \right)^*.
\]
Thus by Proposition~\ref{prop:Nak*}, $I^{[q]} \subseteq (J^{[q]})^*$,
which implies that $I \subseteq J^*$ since $R$ has a weak test element.
In particular, $z_1 \in J^*$, contradicting the fact that $z_1, \dotsc, z_n$ are $*$-independent.
Hence, $z \in \m I$, so $\sptc{I} \cap I \subseteq \m I$.
\end{proof}

Next, note that $*$-independence behaves well with
respect to Frobenius powers and local module-finite extensions:

\begin{lemma}\label{lem:*indext}
Let $(R,\m)$ be a Noetherian local ring of characteristic $p>0$ containing a $q_0$-weak test
element $c$, let $q'$ be any power of $p$, and let $(S,\n)$ be a local module-finite
extension ring of $R$.  Let $f_1, \dotsc, f_n$ be $*$-independent elements of $R$.  Then $f_1^{q'},
\dotsc, f_n^{q'}$ is a $*$-independent set of elements of $S$.
\end{lemma}

\begin{proof}
Suppose that $f_i^{q'} \in ((f_1^{q'}, \dotsc, \hat{f_i^{q'}}, \dotsc, f_n^{q'}) S)^*$.  Then for all $q \gg 0$,
$c f_i^q \in (f_1^q, \dotsc, \hat{f_i^q}, \dotsc, f_n^q) S$, which implies that \[
f_i \in ((f_1, \dotsc, \hat{f_i}, \dotsc, f_n)S)^* \cap R = (f_1, \dotsc, \hat{f_i}, \dotsc, f_n)^*,
\]
contradicting the $*$-independence of $f_1, \dotsc, f_n$ in $R$.
\end{proof}

We will also make crucial use of the following
\begin{lemma}\label{lem:sp*contain}
If $R$ is a Noetherian local ring of characteristic $p>0$ with a weak test element $c$, then
for all ideals $I \subseteq R$, \[
\sptc{(I^*)} = \sptc{I} = (\sptc{I})^*.
\]
Hence, if $J^* \subseteq I^*$, it follows that $\sptc{J} \subseteq \sptc{I}$.
\end{lemma}

In other words, $\sptc{I}$ is tightly closed, and it coincides with $\sptc{J}$ whenever
$I^* = J^*$.

\begin{proof}
Let $f \in \sptc{(I^*)}$.  Then there is some $q_0$ such that for all $q \gg 1$, \[
c f^q \in \m^{[q / q_0]} (I^*)^{[q]}.
\]
Multiplying both sides by $c$, we get: \[
c^2 f^q \in \m^{[q / q_0]} (c (I^*)^{[q]}) \subseteq \m^{[q / q_0]} I^{[q]}.
\]
Hence $f^{q_0} \in (\m I^{[q_0]})^*$, so that $f \in \sptc{I}$.

Now suppose $f \in (\sptc{I})^*$.  Then for $q \gg 0$, $c f^q \in (\sptc{I})^{[q]}$.  However,
there is some $q_0$ such that for $q \gg q_0$, $c (\sptc{I})^{[q]} \subseteq \m^{[q/q_0]} I^{[q]}$.
Thus, we have $c^2 f^q \in c (\sptc{I})^{[q]} \subseteq \m^{[q/q_0]} I^{[q]}$, which
means that $f \in \sptc{I}$.

For the last statement, we have $\sptc{J} = \sptc{(J^*)} \subseteq \sptc{(I^*)} = \sptc{I}$.
\end{proof}

\section{Special tight closure decomposition}\label{sec:sptcx}
In this section, we extend \cite[Theorem 2.1]{HuVr} of Huneke and Vraciu.
First, recall two lemmas from their paper:

\begin{lemma}\label{lem:HuVr2.3}\cite[Lemma 2.3]{HuVr}
Let $(R,\m)$ be an excellent normal local ring of characteristic $p>0$
and let $I = (f_1, \dotsc, f_n)$ be an ideal.  If $f \in I^*$, there is
a test element $c$ and a power $q_0$ of $p$ such that \[
c f^q \in c I^{[q]} + \m^{[q / q_0]} I^{[q]}
\]
for all $q \geq q_0$. \footnote{In \cite{HuVr}, this lemma is stated
without the brackets on the power of $\m$.  However, as the reader
can easily verify, their statement and the present one are equivalent, albeit
with a different choice of $q_0$.}
\end{lemma}

\begin{lemma}\label{lem:HuVr2.2}\cite[Lemma 2.2]{HuVr}
Let $(R,\m)$ be an excellent analytically irreducible local ring
of characteristic $p>0$.  Let $I = (f_1, \dotsc, f_n)$ and $f \in I^*$.
Assume that for each
$i = 1, \dotsc, n$, there exists $\alpha_i \in R$ such that \[
f_i \notin (f + \alpha_i f_i, f_1, \dotsc, f_{i-1}, f_{i+1}, \dotsc, f_n)^*.
\]
Then $f \in I + \sptc{I}$.
\end{lemma}

We will also use the following two results.  The first of these is known, but
I can't find the explicit statement anywhere in the literature.  Thanks to William Heinzer
for pointing out its proof:

\begin{lemma}\label{lem:normlocal}
Let $R$ be an analytically irreducible Noetherian local ring.  Then $\bar{R}$, its
normalization, is also local.
\end{lemma}

\begin{proof}
According to Nagata \cite[Theorem 43.20]{NagLR}, for any quasi-local integral domain $R$,
there is a one-to-one correspondence between the maximal ideals of $\bar{R}$ and
the associated primes of the Henselization $R^h$ of $R$.  However, in our case,
since $R$ is Noetherian, $R^h$ is a subring of $\hat{R}$, which is an integral
domain by assumption.  Thus $R^h$ is also an integral domain, so it has only one
associated prime, which then implies that $\bar{R}$ has only one maximal ideal.
\end{proof}

\begin{notation} In the next lemma and for the rest of this paper, if $R$ is a local
ring, $\kappa(R)$ denotes its residue field.
\end{notation}

\begin{lemma}\label{lem:nIcapR}
Let $(R,\m) \subseteq (S,\n)$ be a local module-finite extension of Noetherian local
rings of characteristic $p>0$.  Then for any ideal $I \subseteq R$,
$\sptc{(I S)} \cap R = \sptc{I}$.  Also, if $\kappa(R) = \kappa(S)$, then
$I S \cap R \subseteq I + \sptc{I}$.
\end{lemma}

\begin{proof}
First of all, since the special part of the tight closure can be computed modulo minimal primes,
we may assume that $R$ and $S$ are reduced.  Then by Hochster \& Huneke \cite[Proposition 5.19]{HHsplit},
there is an $R$-linear map $\phi: S \rightarrow R$ such that $d = \phi(1) \in R^o$.
Let $q_1$ be a power of $p$ such that $\n^{[q_1]} \subseteq \m S$.  Suppose
$f \in \sptc{(I S)} \cap R$.  There is some $b \in S^o$ and some power $q_0$ of $p$ such
that for all powers $q \gg 0$ of $p$, $b f^q \in \n^{[q/q_0]} I^{[q]}$.  Moreover,
by the second paragraph of the proof of \cite[Corollary 5.23]{HHsplit}, there is some
$s \in S$ such that $s b = c \in R^o$.

Then \[
d c f^q = \phi(c f^q) \in \phi(\n^{[q/q_0]} I^{[q]}) \subseteq
 \phi(\m^{\left[\frac{q} {q_0 q_1}\right]} I^{[q]} S) \subseteq \m^{\left[\frac{q} {q_0 q_1}\right]} I^{[q]},
\]
which implies that $f^{q_0 q_1} \in (\m I^{[q_0 q_1]})^*$, whence $f \in \sptc{I}$.

For the second statement, let $f \in I S \cap R$.  Say $I = (f_1, \dotsc, f_r)$.
Then $f = \sum_{i=1}^n s_i f_i$ for some $s_i \in S$.  If $\kappa(S) = \kappa(R)$,
then for each $i$ there exists $r_i \in R$ and $\mu_i \in \n$ such that
$s_i = r_i + \mu_i$.  Then \[
\sum_{i=1}^n \mu_i f_i = f - \sum_{i=1}^n r_i f_i \in \n I \cap R \subseteq \sptc{I}
\]
by the first part of the lemma.  Hence, since $\sum_{i=1}^n r_i f_i \in I$,
we have $f \in I + \sptc{I}$.
\end{proof}

Now we come to the promised extension of Huneke and Vraciu's theorem.
They prove the result in the case where $R$ is normal and $k$ is perfect.  In the
following, for an ideal $J$ of $R$, the notation $J^{1/q}$ will denote the ideal of elements of
$R$ whose $q$'th powers are in $J$.

\begin{thm}\label{thm:huvrq}
Let $(R,\m,k)$ be an excellent analytically irreducible local domain of
characteristic $p>0$.  Assume that $k = \kappa(\bar{R})$.  Then for any proper
ideal $I$ of $R$, there is some $q'$ such that \begin{equation}\label{eq:huvrq}
I^* = \left(I^{[q']} + \sptc{(I^{[q']})}\right)^{1/q'}.
\end{equation}
If either (a) $k$ is perfect or (b) $\dim R \leq 1$ and $k$ is infinite,
we may choose $q' = 1$.  In case (b), even if $k \neq \kappa(\bar{R})$ we
have $I^* = I \bar{R} \cap R$.
\end{thm}

\begin{proof}
The containment `$\supseteq$' follows from the definitions:  Let
$x \in \left(I^{[q']} + \sptc{(I^{[q']})}\right)^{1/q'}$ for some $q'$.  Then
$x^{q'} = y + z$, $y \in I^{[q']}$, $z \in \sptc{(I^{[q']})}$, so \[
c x^q = c y^{q/q'} + c z^{q/q'} \subseteq I^{[q]} + \m^{[q/q_0]} I^{[q]} = I^{[q]}
\]
for some $q_0$ and sufficiently high $q$, which implies that $x \in I^*$.

Moreover, if $J$ is a $*$-reduction
of $I$ and the theorem holds for $J$, then \[
I^* = J^* \subseteq \left(J^{[q']} + \sptc{(J^{[q']})}\right)^{1/q'}
\subseteq \left(I^{[q']} + \sptc{(I^{[q']})}\right)^{1/q'}.
\]
Hence we may assume without loss of generality that $I$ is generated by a $*$-independent set.
Also, if $\dim R = 0$ then $R$ is a field.  So we may also assume $\dim R \geq 1$.

First consider the case where $\dim R = 1$ and $k$ is infinite.  In this case, $I$ has a principal reduction $(g)$, since the analytic spread of an
ideal is bounded by the dimension of its ring.
Hence, $I \subseteq \intcl{(g)} = (g)^*$, since tight and integral closure coincide
for principal ideals.  Thus we may
assume that $I = (g)$.  But then $\intcl{I} = \intcl{(g)} = g \bar{R} \cap R$.
This proves the last statement of the theorem, and if $k = \kappa(\bar{R})$, then
Lemma~\ref{lem:nIcapR} yields Equation~(\ref{eq:huvrq}) with $q'=1$, in case (b).

Now we treat case (a) as well as the general case.
Let $f \in I^*$, and let $f_1, \dotsc, f_n$ be a $*$-independent generating
set of $I$.  Then by Lemma~\ref{lem:HuVr2.3}, for all $q \geq q_0$
there exists $u_q \in \bar{R}$ such that \begin{equation}\label{eq:first}
c f^q - c u_q f_1^q \in (f_2, \dotsc, f_n)^{[q]} \bar{R} + \n^{[q / q_0]} f_1^q,
\end{equation}
where $\n$ denotes the maximal ideal of $\bar{R}$.  Taking this equation to the $p$'th power:
\[
c^p \left(f^{p q} - (u_q)^p f_1^{p q}\right) \in (f_2, \dotsc, f_n)^{[p q]} \bar{R}
 + \n^{[p q / q_0]} f_1^{p q}
\]
On the other hand, replacing $q$ by $p q$ in (\ref{eq:first}) and multiplying by $c^{p-1}$:
\[
c^p \left(f^{p q} - u_{p q} f_1^{p q}\right) \in (f_2, \dotsc, f_n)^{[p q]} \bar{R}
 + \n^{[p q / q_0]} f_1^{p q}.
\]
Subtracting the latest two displayed equations from each other, we have that for each $q$, there exists
$m_q \in \n^{[p q / q_0]}$ such that\[
c^p (u_{p q} - (u_q)^p) - m_q \in (f_2, \dotsc, f_n)^{[p q]} \bar{R} : f_1^{p q} \subseteq \n^{[p q / q_0]}
\]
by Proposition~\ref{prop::} (perhaps requiring an increase in $q_0$), since the $f_i^{p q}$ are
$*$-independent in $\bar{R}$ by Lemma~\ref{lem:*indext}.  (Proposition~\ref{prop::} applies because
$\bar{R}$ is a normal domain, thus analytically irreducible.)  Hence, \[
u_{p q} - (u_q)^p \in \n^{[p q / q_0]} : c^p.
\]
For large enough $q$, this colon ideal is contained in $\n$, by the Artin-Rees lemma
and the fact that $c^p$ is a nonzerodivisor in $\bar{R}$.  Thus
there is some power $q_1$ of $p$ such that,
$u_{q} - (u_{q_1})^{q/q_1} \in \n$ for all powers $q \geq q_1$ of $p$.  Since $k = \kappa(\bar{R})$, there is some $\alpha \in R$ such that
$\alpha + u_{q_1} \in \n$.  Combining these two facts together, we get \[
u_{q} + \alpha^{q / q_1} \in \n \text{ for all } q \geq q_1.
\]
If $k$ is perfect, then modulo $\m$, $\alpha$ has a $q_1$'th root $\alpha'$.  That is, $\alpha - (\alpha')^{q_1} \in \m$.
In this case, we can replace $q_1$ by $1$ and $\alpha$ by $\alpha'$ in the rest of the proof.

Now \emph{suppose that } $f_1^{q_1} \in \left((f^{q_1} + \alpha f_1^{q_1}, (f_2, \dotsc, f_n)^{[q_1]}) \bar{R}\right)^*$.

Then by Lemma~\ref{lem:HuVr2.3} again, \[
c f_1^q \in ((f_2, \dotsc, f_n)^{[q]}\bar{R}, c(f^q + \alpha^{q/q_1} f_1^q)\bar{R}, \n^{[q/q_0]} f^q, \n^{[q/q_0]} f_1^q).
\]
Multiplying by $c$ and using (\ref{eq:first}) we have \[
c^2 f_1^q \in ((f_2, \dotsc, f_n)^{[q]} \bar{R}, c^2(u_q + \alpha^{q/q_1})f_1^q \bar{R}, n^{[q/q_0]} f_1^q).
\]
That is, there exist $B_q \in \bar{R}$ and $C_q \in \n^{[q/q_0]}$ such that \[
(c^2(1 - B_q(u_q + \alpha^{q/q_1})) - C_q) f_1^q \in (f_2, \dotsc, f_n)^{[q]} \bar{R}.
\]
Hence, using Proposition~\ref{prop::} again (along with the fact that the $f_i$ are $*$-independent in $\bar{R}$),
$1 - B_q(u_q + \alpha^{q/q_1}) \in \n^{[q/ q_0]} : c^2$
for all large $q$.  But this latter colon must be contained in $\n$ for large $q$, and then
since $u_q + \alpha^{q/q_1} \in \n$ as well, we arrive at the contradiction that $1\in \n$.

We have shown by contradiction that \[
f_1^{q_1} \notin \left((f^{q_1} + \alpha f_1^{q_1}, (f_2, \dotsc, f_n)^{[q_1]}) \bar{R}\right)^*.
\]
In particular, $f_1^{q_1} \notin (f^{q_1} + \alpha f_1^{q_1}, (f_2, \dotsc, f_n)^{[q_1]})^*$.
Let $\alpha_1 = \alpha$.  By symmetry, for each $1 \leq i \leq n$,  a power $q_i$ of $p$ (which, if
$k$ is perfect, can be chosen as before to be $1$)
and $\alpha_i \in R$ such that \[
f_i^{q_i} \notin \left(f^{q_i} + \alpha_i f_i^{q_i}, (f_1, \dotsc, \hat{f_i}, \dotsc, f_n)^{[q_i]}\right)^*.
\]
Now let $q = \max \{ q_i \mid 1 \leq i \leq n\}$, and let $\beta_i = \alpha_i^{q / q_i}$.  Then we have \[
f_i^q \notin \left( f^q + \beta_i f_i^q, (f_1, \dotsc, \hat{f_i}, \dotsc, f_n)^{[q]}\right)^*
\]
for all $i$.  Then by Lemma~\ref{lem:HuVr2.2}, $f^{q} \in I^{[q]} + \sptc{(I^{[q]})}$,
which means that $f \in \left(I^{[q]} + \sptc{(I^{[q]})}\right)^{1/q}$.

If $k$ is perfect, then $q=1$ in the above, so we are done.  Otherwise, to finish
the proof, note that what we have shown is that \[
I^* \subseteq \bigcup_q \left(I^{[q]} + \sptc{(I^{[q]})}\right)^{1/q}.
\]
However, the ideals in the union are all nested by inclusion, so by the ascending chain condition,
there is some $q'$ for which the $q'$'th ideal equals the $q$'th ideal for all $q \geq q'$, and
this common ideal contains the rest of them.  Hence, \[
I^* \subseteq \left(I^{[q']} + \sptc{(I^{[q']})}\right)^{1/q'},
\]
as was to be shown.
\end{proof}

We now show that the theorem is sharp.  First of all, if $k \neq \kappa(\bar{R})$ and
$k$ is perfect, there
will be a principal ideal $I$ of $R$ such that $I^* \neq I + \sptc{I}$.  To see
this, let $\alpha \in \kappa(\bar{R}) \setminus k$.  Then there are $f, g \in R$ with
$f \in \intcl{(g)}$ such that the image of $f/g$ in $\kappa(\bar{R})$ is $\alpha$.  Suppose
that $f \in (g) + \sptc{(g)}$.  Then there is some
$r \in R$ such that $f - r g \in \sptc{(g)}$,
and there is some $q_0$ such
that for all $q \geq q_0$ \[
c (f - r g)^q \in \m^{[q/q_0]} g^{q} \subseteq \n^{[q/q_0]} g^{q}
\]
Then in $\bar{R}$, we have $c (\frac{f}{g} - r)^q \in \n^{[q/q_0]}$.  Thus,
$(\frac{f}{g} - r)^{q_0} \in \n^* = \n$, so that $\frac{f}{g} - r \in \n$.  In other words,
the image $\alpha$ of $\frac{f}{g}$ in $\kappa(\bar{R})$ lies in $k$, contradicting
our choice of $\alpha$.

Next, note that if $k$ is not perfect, $I^*$ may not equal $I + \sptc{I}$, even if $R$
is normal, as we show in the following example in dimension 2:

Let $T = \frac{\Z}{(p)}[U,V,W,X,Y,Z]$, where $p$ is a positive prime number, and let
$S = T / (U X^p + V Y^p + W Z^p)$.  We use the lower-case letters
$u, \dotsc, z$ to denote the images of $U, \dotsc, Z$ in $S$.  Then $S$ is a
general grade reduction of $(k[X,Y,Z], (X^p,Y^p,Z^p))$ in the sense of Hochster \cite{HoGGR},
so by \cite[Theorem (c)]{HoGGR}, $S$ is an integrally closed domain.  Now
let \[
R = S_{(x,y,z)} \cong \frac{\frac{\Z}{(p)}(u,v,w)[X,Y,Z]_{(X,Y,Z)}} {(X^p + (v/u) Y^p + (w/u) Z^p)}.
\]
Note that as the localization of an
integrally closed domain at a prime ideal, $R$ is normal, with maximal ideal $\m = (x,y,z)$ and
residue field $k = \frac{\Z}{(p)}(u,v,w)$.  Let $I = (y,z) \subseteq R$.
Then \[
x^p = -(v/u) y^p - (w/u) z^p \in I^{[p]},
\] so that $x \in I^F \subseteq I^*$.  However,
I claim that $x \notin I + \sptc{I}$.  For suppose that $x \in I + \sptc{I}$.  Then
there exist $r,s \in R$ such that \[
r y + s z - x \in \sptc{I}.
\]  In particular, there exists
some $q_0$ such that for $q \gg q_0$, there exist $m_q, n_q \in \m^{[q/q_0]}$ such that \begin{align*}
c \left( \left(\frac{v}{u}\right)^{q/p} + r^q\right) y^q + c \left( \left(\frac{w}{u}\right)^{q/p} + s^q\right) z^q
 &= c(r^q y^q + s^q z^q - x^q) \\
 &= c (r y + s z - x)^q \\
 &= m_q y^q + n_q z^q.
\end{align*}
Then since $y,z$ are $*$-independent (since they form a system of parameters for $R$), the colon
criterion shows that there is some $q_1 \geq q_0$
such that for all $q \geq q_1$, \[
c \left(\left(\frac{v}{u} + r^p\right)^{q_1/p}\right)^{q/q_1} = c \left(\left(\frac{v}{u}\right)^{q/p} + r^q\right) \in \m^{[q/q_1]}.
\]
Therefore, $(v/u + r^p)^{q_1/p} \in \m^* = \m$, so that since $\m$ is radical, $v/u + r^p \in \m$.
Hence, the image of $r$ in $R / \m \cong k$ is a $p$'th root of $-v/u$, which is a contradiction
since $-v/u$ has no $p$'th root in $k$.

\section{$*$-spread}\label{sec:*spread}

The following is the main theorem of this note.  It extends a result of Vraciu
\cite[Proposition 3.6]{Vr*ind} which covers the case where both $I / J$ and $I / K$ have finite length.

\begin{thm}\label{thm:*spread}
Let $(R,\m,k)$ be an excellent analytically irreducible local domain of characteristic $p>0$.
Let $I$ be an ideal of $R$, and let $J$, $K$ be minimal $*$-reductions of $I$.  Then
$\mu(J) = \mu(K)$.
\end{thm}

\begin{proof}
Without loss of generality, $\mu(J) \leq \mu(K)$.  Let $a_1, \dotsc, a_n$ be a
minimal generating set for $J$, and let $b_1, \dotsc, b_r$ be a minimal generating set
for $K$. So $r \geq n$, and we want to prove equality.

We will prove the following by induction on $n - i$:

\noindent \textbf{Claim: } We can reorder the $a_j$'s in such a way that
whenever $0 \leq i \leq n$, $I \subseteq (a_1, \dotsc, a_i, b_{i+1}, \dotsc, b_n)^*$.\footnote{
Note that this method of proof shows that in fact, $I$ is a \emph{matroid}
(a notion familiar to combinatorialists) whose independent
sets are the $*$-independent subsets of $I$.  Joe Brennan suggested this interpretation when
I told him about $*$-spread.}

\begin{proof}[Proof of Claim.]
For short, let $L_i = (a_1, \dotsc, a_i, b_{i+1}, \dotsc, b_n)$.

If $i=n$, there is nothing to prove.  So assume $1 \leq i \leq n$ and the claim holds
for $i$.  It is enough to show that for some reordering of
$a_1$ through $a_i$,
$a_i \in L_{i-1}^*$ for if we can show that, then we have \[
I \subseteq L_i^* = ((a_1, \dotsc, a_{i-1}, b_{i+1}, \dotsc, b_n) + (a_i))^*
\subseteq (L_{i-1}^* + (a_i))^* = L_{i-1}^*.
\]

We have for some $q_1$ that \[
b_i \in I^* = L_i^* = (L_i \bar{R})^* \cap R = \left((L_i^{[q_1]} \bar{R}) + \sptc{(L_i^{[q_1]} \bar{R})}\right)^{1/q_1} \cap R,
\]
where the equality is by Theorem~\ref{thm:huvrq}, so that \[
c b_i^q = \sum_{j=1}^i (c r_j^{q/q_1} + t_{j,q}) a_j^q + \sum_{k=i+1}^n (c d_k^{q/q_1} + u_{k,q}) b_k^q,
\]
where $r_j, d_k \in \bar{R}$ and $t_{j,q}, u_{k,q} \in \n^{[q/q_0]}$ for some fixed power
$q_0$ of $p$, for all $q \geq \max\{q_0,q_1\}$.
Here $\n$ is the maximal ideal of $\bar{R}$.  Also we can replace $q_0$ and $q_1$ by $\max\{q_0,q_1\}$
and assume without loss of generality that $q_0 = q_1$.

Now suppose that $r_j \in \n$ for all $1 \leq j \leq i$.  Then \[
c (b_i^{q_0} - \sum_{k=i+1}^n d_k b_k^{q_0})^{q/q_0} = c b_i^q - \sum_{k=i+1}^n c d_k^{q/q_0} b_k^q \in
\n^{[q/q_0]} L_i^{[q]} = (\n L_i^{[q_0]})^{[q/q_0]}.
\]
Hence, $b_i^{q_0} - \sum_{k=i+1}^n d_k b_k^{q_0} \in (\n L_i^{[q_0]})^* \subseteq \sptc{(L_i^{[q_0]} \bar{R})}$.  Moreover,
we have that $\sptc{(L_i^{[q_0]} \bar{R})} = \sptc{(K^{[q_0]} \bar{R})}$ since $L_i^* = K^*$.  This
follows from Lemma~\ref{lem:sp*contain} in addition to the fact
that tight closure is persistent.  Therefore: \begin{align*}
b_i^{q_0} - \sum_{k=i+1}^n d_k b_k^{q_0} &\in \sptc{(L_i^{[q_0]} \bar{R})} \cap K^{[q_0]} \bar{R} \\
&= \sptc{(K^{[q_0]} \bar{R})} \cap K^{[q_0]} \bar{R} \\
&= \text{ (by Lemma~\ref{lem:sp*capI}) } \n K^{[q_0]} \subseteq \n I^{[q_0]}.
\end{align*}
But this contradicts the fact that $b_i^{q_0}, \dotsc, b_n^{q_0}$ are linearly independent
modulo $\n I^{[q_0]}$ (which is true since $b_1^{q_0}, \dotsc, b_n^{q_0}$ are $*$-independent as
elements of $\bar{R}$.)

Hence there is some $j$ with $1 \leq j \leq i$ and $r_j \notin \n$.  By rearranging
the first $i$ indices, we may assume that $r_i \notin \n$.  Then we have \begin{align*}
ca_i^q &= r_i^{-q/q_0} \left( c b_i^q - \sum_{k=i+1}^n (c d_k^{q/q_0} - u_{k,q}) b_k^q
 - \sum_{j=1}^{i-1} (c r_j^{q/q_0} + t_{j,q}) a_j^q \right) \\
 & \quad \quad - r_i^{-q/q_0} t_{i,q} a_i^q &\\
&\in (L_{i-1} \bar{R})^{[q]} + \n^{[q / q_0]} a_i^q.& &
\end{align*}

Then by the Nakayama lemma for tight closure, we have \[
a_i \in (L_{i-1} \bar{R})^* \cap R = L_{i-1}^*,
\]
as required.
\end{proof}

By the Claim, taken with $i=0$, $I \subseteq L_0^*$, so that $L_0 = (b_1, \dotsc, b_n)$
is a $*$-reduction of $I$.  But since $K$ is minimal with respect to this property and $L_0 \subseteq K$, we have
$K = L_0$, which forces $n=r$.
\end{proof}

\section*{acknowledgements}
The author wishes to thank his thesis advisor, Craig Huneke, for many
useful conversations and constant support, and also Dan Katz for teaching the course that inspired
this work.

\bibliographystyle{amsalpha}
\providecommand{\bysame}{\leavevmode\hbox to3em{\hrulefill}\thinspace}
\providecommand{\MR}{\relax\ifhmode\unskip\space\fi MR }
\providecommand{\MRhref}[2]{%
  \href{http://www.ams.org/mathscinet-getitem?mr=#1}{#2}
}
\providecommand{\href}[2]{#2}

\end{document}